\documentclass[10pt,fleqn]{article}

\usepackage{ifpdf}
\usepackage{amsmath,amssymb,amsthm}

\ifpdf
    \usepackage[pdftex]{graphicx}
    \usepackage[colorlinks=true, pdfstartview=FitH, linkcolor=blue,
        citecolor=blue, urlcolor=blue]{hyperref}
\else
   \usepackage{graphicx}
   \usepackage[dvips]{hyperref}
\fi
\usepackage[numeric]{amsrefs} 

    \newcommand{\MathReview}[1]{~\href{http://www.ams.org/mathscinet-getitem?mr=#1}{\mbox{\bf MR~#1}}}

    \newcommand{\NN}{\ensuremath{\mathbb N}}
    \newcommand{\ZZ}{\ensuremath{\mathbb Z}}

    \newcommand{\Spec}[1]{\text{\sc Spec\ensuremath{(#1)}}}

    \newcommand{\floor}[1]{\mbox{$\left\lfloor #1 \right\rfloor$}}
    
    \newcommand{\fp}[1]{\mbox{$\left\{ #1 \right\}$}}

    \newtheorem{thm}{Theorem}[section]
    \newtheorem{lem}[thm]{Lemma}
    
    \newtheorem{cor}[thm]{Corollary}

\title{Gaps in the Spectrum of Heights of Projective Points}
\author{
    Kevin O'Bryant\thanks{Supported by PSC-CUNY grant 60070-36 37}\\
    {\small City University of New York, College of Staten Island, New York, NY}\\
    {\small \tt{kevin@member.ams.org}}}
\date{\today}

\begin{document}
    \ifpdf
       \DeclareGraphicsExtensions{.pdf,.jpg,.mps,.png}
    \fi
\maketitle \sloppypar

\begin{abstract}
Let $\ZZ_m$ be the ring of integers modulo $m$ (not necessarily prime), $\ZZ_m^\ast$ its multiplicative group, and let $x\bmod m$ be the least nonnegative residue of $x$ modulo $m$. The {\em height} of a point $r=\langle r_1,\dots,r_d\rangle\in (\ZZ_m^\ast)^d$ is $h_m(r)=\min\left\{ \sum_{i=1}^d (k r_i \bmod m) \colon k=1,\dots,p-1\right\}$. For $d=2$, we give an explicit formula in terms of the convergents to the continued fraction expansion of $\bar r_1 r_2 /m$. Further, we show that the multiset $\{m^{-1} h_m((r_1,r_2)): m\in\NN, r_i\in\ZZ_m^\ast \}$, which is trivially a subset of $[0,2]$, has only the numbers $1/k$ ($k\in \ZZ^+$) and 0 as accumulation points.
\end{abstract}

\section{Introduction}
In \cite{Nathanson.Blair}, Nathanson \& Sullivan raised the problem of bounding the height of points in $(\ZZ_m^\ast)^d$, where $m$ is a prime. After proving some general bounds for $d>2$, they move to identifying those primes $p$ and residues $r$ with $h_p(\langle 1,r \rangle)>(p-1)/2$. In particular, they prove that if $h_p(\langle 1,r\rangle )<p$, then it is in fact at most $(p+1)/2$. Nathanson has further proven \cite{Nathanson} that if $p$ is a sufficiently large prime and $h_p(\langle 1,r\rangle)<(p+1)/2$, then it is in fact at most $(p+4)/3$. In other words, $p^{-1}h_p(\langle 1,r\rangle)$ is either near 1, near $1/2$, or at most $1/3$.

In this paper we show that these gaps in the values of $p^{-1}h_p(\langle 1,r\rangle)$ continue all the way to 0, even if $p$ is not restricted to be prime. The main tool is the simple continued fraction of $r/p$. To avoid confusion, as we do not use primeness here, and since the numerators of continued fractions are traditionally denoted by $p$, we denote our modulus by $m$.

\begin{figure}
\begin{picture}(390,250)
    \put(20,15){\includegraphics[width=5in]{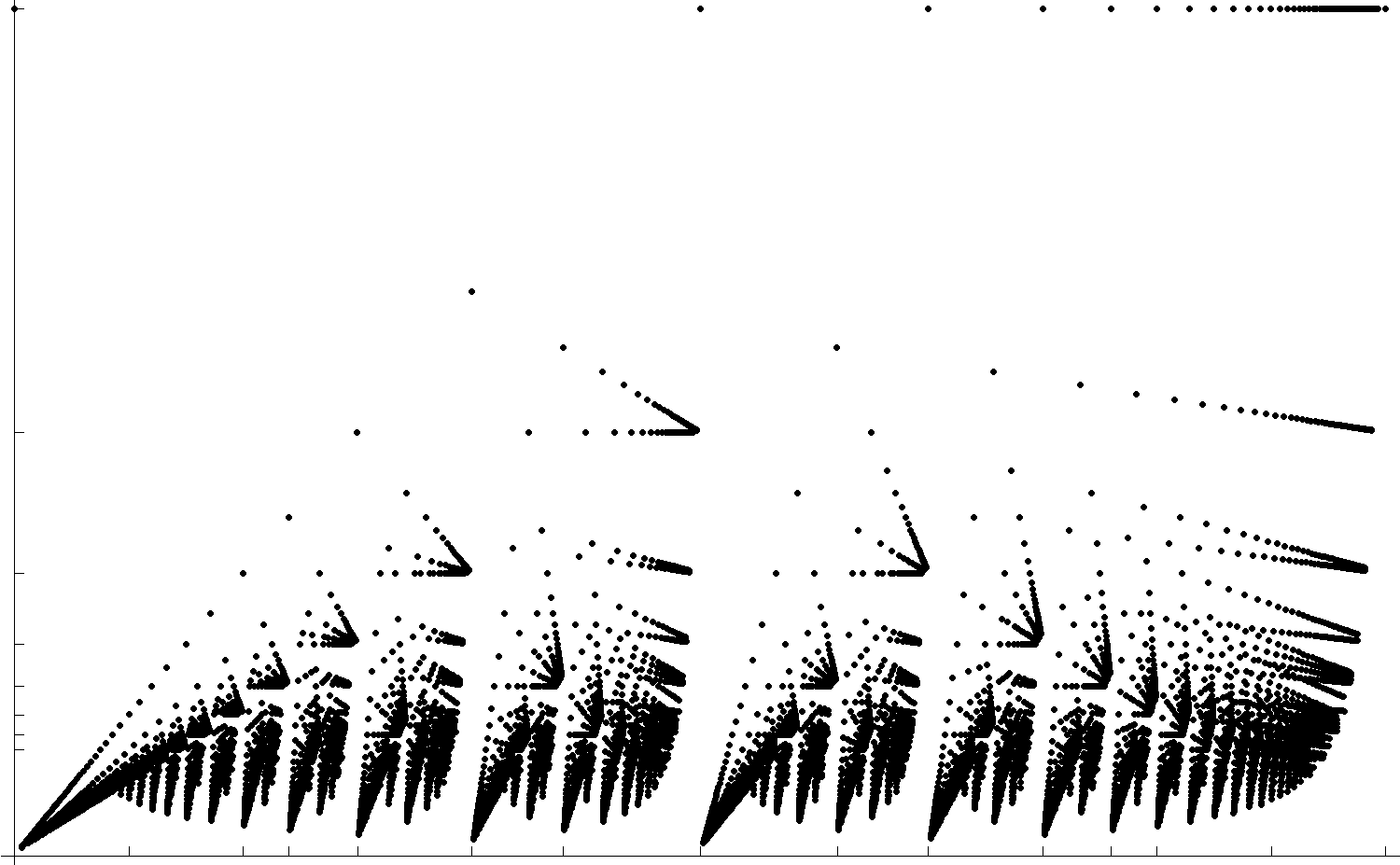}}
    \put(10,15){0}
    \put(0,41){$1/8$}
    \put(0,51){$1/6$}
    \put(0,59){$1/5$}
    \put(0,69){$1/4$}
    \put(0,87){$1/3$}
    \put(0,123){$1/2$}
    \put(10,232){1}
    \put(21,0){0}
    \put(48,0){$\tfrac 1{12}$}
    \put(79,0){$\tfrac 1{6}$}
    \put(91,0){$\tfrac 1{5}$}
    \put(110,0){$\tfrac 1{4}$}
    \put(138,0){$\tfrac 1{3}$}
    \put(162,0){$\tfrac 2{5}$}
    \put(198,0){$\tfrac 1{2}$}
    \put(233,0){$\tfrac 3{5}$}
    \put(257,0){$\tfrac 2{3}$}
    \put(287,0){$\tfrac 3{4}$}
    \put(304,0){$\tfrac 4{5}$}
    \put(316,0){$\tfrac 5{6}$}
    \put(344,0){$\tfrac {11}{12}$}
    \put(375,0){1}
\end{picture}
\caption{The points $(\frac rm, H(\frac rm))$, for all $0<r<m\le 200$.\label{Fareypic}}
\end{figure}
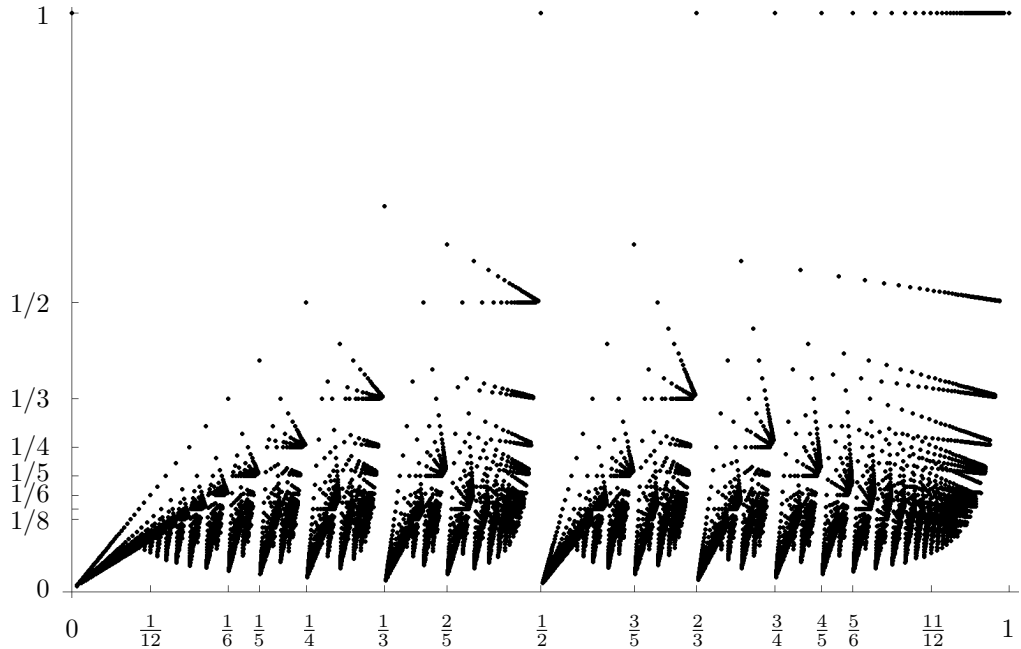

If $\gcd(r_1,m)=1$, then $h_m(\langle r_1,r_2\rangle)=h_m(\langle 1,\bar r_1\,r_2\rangle)$, and so we may assume without loss of generality that $r_1=1$. We are thus justified in making the following definition for relatively prime positive integers $r,m$:
    \begin{align*}
    H(r/m) &:= m^{-1} \cdot h_m(\langle 1,r \rangle) \\
        &= m^{-1} \cdot \min\{ k + (kr \bmod m) \colon 1\le k <m\}\\
        &=\min \{ k/m + \fp{kr/m} \colon 1\le k < m\}.
    \end{align*}
We are using the common notation for the fractional part of $x$, namely $\fp{x}:=x-\floor{x}$. Figure~\ref{Fareypic} shows the points $(\frac rm, H(\frac rm))$ for all $r,m\le 200$.

The spectrum of a set $M\subseteq\NN$, written \Spec{M}, is the set of real numbers $\beta$ with the property that there are $m_i\in M$, $m_i\to\infty$, and a sequence $r_i$ with $\gcd(r_i,m_i)=1$, and $H(r_i/m_i)\to \beta$. Nathanson~\cite{Nathanson} and Nathanson~\&~Sullivan~\cite{Nathanson.Blair} proved that
    \[\Spec{\text{primes}}\cap [\tfrac 13,\infty) = \left\{\tfrac 13, \tfrac 12, 1\right\}.\]

Our main theorem concerns the spectrum of heights.
\begin{thm}\label{MainTheorem}
$\Spec{\NN} =\{0\}\cup\{1/k \colon k\in\ZZ^+\}$.
\end{thm}

\section{Continued Fractions}
For a rational number $0<\frac rm<1$, let $[0;a_1,a_2, \dots,a_n]$ be (either one of) its simple continued fraction expansion, and let $p_k/q_k$ be the $k$-th convergent. In particular
    \begin{align*}
    \frac{p_0}{q_0} & = \frac 01 \\
    \frac{p_2}{q_2} & = \frac {a_2}{1+a_1a_2} \\
    \frac{p_4}{q_4} & = \frac {a_2+a_4+a_2 a_3 a_4}{1+a_1 a_2+a_1 a_4+a_3 a_4+a_1 a_2 a_3 a_4}
    \end{align*}
The $q_i$ satisfy the recurrence $q_{-2}=1,q_{-1}=0,q_n=a_n q_{n-1}+q_{n-2}$ (with $a_0=0$), and are called the {\em continuants}. The {\em intermediants} are the numbers $\alpha q_{n-1}+q_{n-2}$, where $\alpha$ is an integer with $1\le \alpha\leq a_n$.

Let $E[a_0,a_1,\dots,a_n]$ be the denominator $[a_0;a_1,\dots,a_n]$, considered as a polynomial in $a_0,\dots,a_n$, and set $E[]=1$. Then $p_k=E[a_0,\dots,a_k]$ and $q_k=E[a_1,\dots,a_k]$. We will make use of the following combinatorial identities, which are in \cite{Chapter 13}*{Roberts}, with $0<s<t<n$:
    \[q_\ell=q_k E[a_{k+1},\dots,a_\ell]+q_{k-1}E[a_{k+2},\dots,a_\ell],\]
    \[p_n E[a_s,\dots,a_t]-p_t E[a_s,\dots,a_n]
        =(-1)^{t-s+1}E[a_0,\dots a_{s-2}]E[a_{t+2},\dots,a_n].\]

The following lemmas are well known. The first is a special case of the ``best approximations theorem'' \cite{HW}*{Theorems 154 and 182}, and the second is an application of~\cite{HW}*{Theorem 150}, the identity $p_nq_{n-1}-p_{n-1}q_n = (-1)^{n-1}$. The third and fourth lemmas follow from the identities for $E$ given above.

\begin{lem}\label{Intermediants}
Fix a real number $x=[0;a_1,a_2,\dots]$, and suppose that the positive integer $\ell$ has the property that $\fp{\ell x} \le \fp{k x}$ for all positive integers $k\le \ell$. Then there are nonnegative integers $n,\alpha\le a_n$ such that $\ell=\alpha q_{2n-1}+q_{2n-2}$.
\end{lem}

\begin{lem}\label{ConsecutiveConvergents}
Let $\frac{p_{2k}}{q_{2k}}=[0;a_1,a_2\dots,a_{2k}]$, and let $x=[0;a_1,a_2\dots,a_{2k-1},a_{2k}+1]$. Then
    \[q_{2k} \cdot x - p_{2k} = \frac{1}{2q_{2k}+q_{2k-1}}.\]
\end{lem}

We will use Fibonacci numbers, although the only property we will make use of is that they tend to infinity): $F_1=1$, $F_2=2$, and $F_n=F_{n-1}+F_{n-2}$.

\begin{lem}\label{ratioofcontinuants}
For all $k\ge1$, $q_k\ge F_k$. For $\ell>k$,
    \[{q_\ell} > {q_k} F_{\ell-k},\qquad \text{ and }\qquad
    q_\ell > a_\ell q_k.\]
\end{lem}

\begin{lem}\label{simplification} For $0<2k<n$,
\[q_{2k} p_n - p_{2k} q_n = E[a_{2k+2},\dots,a_n],\]
where $E[a_{n+1},\dots,a_n]=1$.
\end{lem}

We now state and prove our formula for heights.
\begin{thm}\label{CFbound}
Let $\frac rm=[0;a_1,a_2,\dots,a_n]$ (with $\gcd(r,m)=1$). Then
    \[H(\tfrac rm) = \min_{0\le k < n/2} \;\big\{q_{2k} \tfrac{r+1}{m}-p_{2k}\big\}.\]
\end{thm}

\begin{proof}
First, recall that
    \[H(r/m)= \min\left\{ {k/m}+ \fp{kr/m} \colon 1\le k<m\right\}.\]
Set
    \[I := \{\alpha q_{2i-1}+q_{2i-2} \colon 0\leq \alpha \leq a_{2i}, 0\leq i\leq n/2\}.\]

We call $\ell$ a best multiplier if
    \[{\ell/m}+\fp{\ell r/m} < k/m+\fp{k r/m}\]
for all positive integers $k<\ell$. We begin by proving by induction that the set of best multipliers is contained in the set $I$. Certainly 1 is a best multiplier and also $1=0\cdot q_{-1}+q_{-2}\in I$. Our induction hypothesis is that the best multipliers that are less than $\ell$ are all contained in $I$.

Suppose that $\ell$ is a best multiplier: we know that
    \[\frac{k}m+\fp{k \frac rm} > \frac{\ell}m+\fp{\ell \frac rm}\]
for all $1\le k<\ell$. Since $k<\ell$, we then know that $\fp{k r/m}> (\ell-k)/m+\fp{\ell r/m}>\fp{\ell r/m}$. Lemma~\ref{Intermediants} now tells us that $\ell\in I$. This confirms the induction hypothesis, and establishes that
    \begin{equation}\label{equ:1}
    H(r/m)= \min \{k/m + \fp{kr/m} \colon k\in I\}.
    \end{equation}

Now, note that the function $f_i$ defined by
    \[f_i(x) := \frac{x q_{2i-1}+q_{2i-2}}m + \fp{(x q_{2i-1}+q_{2i-2}) \frac rm}\]
is monotone on the domain $0\le x \le a_{2i}$. As $0 q_{2i-1}+q_{2i-2}=q_{2i-2}$ and $a_{2i}q_{2i-1}+q_{2i-2}=q_{2i}$, this means that the minimum in Eq.~\eqref{equ:1} can only occur at $q_{2i}$, with $0\le 2i \le n$.

As a final observation, we note that $q_0/m+\fp{q_0 r/m} = (r+1)/m$ is at most as large as $q_n/m + \fp{q_n r/m}=1$ (as $q_n=m$). Thus, the minimum in Eq.~\eqref{equ:1} cannot occur exclusively at $k=q_n=m$.
\end{proof}

\begin{cor}\label{Weakbound}
Let $0<r<m$, with $\gcd(r,m)=1$, and let $\frac rm=[0;a_1,\dots,a_n]$, with $a_n\ge 2$. For all $k\in(0,n/2)$,
    \[H(\tfrac rm) \leq \frac{q_{2k}}{m} + \frac{1}{2q_{2k}}.\]
\end{cor}

\begin{proof}
First, note that $\frac rm < [0;a_1,a_2,\dots,a_{2k-1},a_{2k}+1]$. Now, as a matter of algebra (using Lemma~\ref{ConsecutiveConvergents}),
    \begin{multline*}
    q_{2k} \frac{r+1}{m}-p_{2k}
                \le q_{2k}\left( [0;a_1,a_2,\dots,a_{2k}+1]+\frac 1m\right) - p_{2k}
                = \frac{q_{2k}}{m} + \frac{1}{2q_{2k}+q_{2k-1}}\\
                \le \frac{q_{2k}}{m} + \frac{1}{2q_{2k}}.
    \end{multline*}
\end{proof}

\section{Proof of Theorem 1.1}
First, we note that $H(a_2/(1+a_1a_2)) = (1+a_2)/(1+a_1a_2) \to 1/a_1$, where $a_1$ is fixed and $a_2\to\infty$. Thus, $1/k\in \Spec{\NN}$ for every $k$. Also, $H(1/a_1)=2/a_1\to0$ as $a_1\to\infty$, so $0\in\Spec{\NN}$. The remainder of this section is devoted to proving that if $\beta>0$ is in $\Spec{\NN}$, then $\beta$ is rational with numerator 1.

Fix a large integer $s$. Let $r/m$ be a sequence (we will suppress the index) with $\gcd(r,m)=1$ and with $H(r/m)\to \beta> \frac{1}{F_{2s}}$, where $F_{2s}$ is the $2s$-th Fibonacci number: $F_0=0$, $F_1=1$, $F_{i}=F_{i-1}+F_{i-2}$.

Define $a_1,a_2,\dots$ by \[\frac rm = [0;a_1,a_2,\dots,a_n],\] and we again remind the reader that $r/m$ is a sequence, so that each of $a_1, a_2,\dots, $ is a sequence, and $n$ is also a sequence. To ease the psychological burden of considering sequences that might not even be defined for every index, we take this occasion to pass to a subsequence of $r/m$ that has $n$ nondecreasing. Further, we also pass to a subsequence on which each of the sequences $a_i$ is either constant or monotone increasing.

First, we show that $n$ is bounded. Note that $q_{2s}/m$ (fixed $s$) is the same as $q_{2s}/q_n$, and by Lemma~\ref{ratioofcontinuants} this is at most $1/(2F_{2s})$, provided that $n$ is large enough so that $F_{n-2s}>2 F_{2s}$. Take such an $n$. We have from Corollary~\ref{Weakbound} that
    \[
    H(\tfrac rm) \le \frac{q_{2s}}{m} + \frac{1}{2q_{2s}}
                < \frac{1}{2F_{2s}} + \frac{1}{2F_{2s}}
                < \frac{1}{F_{2s}}
                < \beta.
    \]
This contradicts the hypothesis that $H(r/m)\to\beta>0$, and proves that $n$ must be small enough so that $F_{n-2s}>2 F_{2s}$.

Since $m\to\infty$ but $n$ is bounded, some $a_i$ must be unbounded. Let $i$ be the least natural number such that $a_i$ is unbounded.

First, we show that $i$ is not odd. If $i=2k+1$, then
    \[
    H(\tfrac rm) \le q_{2k} \tfrac{r+1}{m}-p_{2k}
    \]
and $p_{2k}$ and $q_{2k}$ are constant. Since $a_{2k+1}\to\infty$, the ratio
    \[\frac rm \to [0;a_1,a_2,\dots,a_{2k}]= \frac{p_{2k}}{q_{2k}}.\]
Thus, since $q_{2k}/m \leq 1/a_{2k+1} \to 0$,
    \[
    H(\tfrac rm) \leq q_{2k} \tfrac{r+1}{m}-p_{2k} = q_{2k} \frac rm + \frac{q_{2k}}{m}-p_{2k}
            \to q_{2k} \frac{p_{2k}}{q_{2k}}+0-p_{2k} = 0,\]
contradicting the hypothesis that $\beta>0$.

Now we show that there are not two $a_i$'s that are unbounded. Suppose that $a_{2k}$ and $a_{j}$ are both unbounded, with $j>2k$. Then
    \[
    H(\tfrac rm) \le \frac{q_{2k}}{m} + \frac{1}{2q_{2k}}.
    \]
Since $a_{2k}$ is unbounded, $\frac{1}{2q_{2k}}\to0$. And since $a_{j}$ is also unbounded,
    \[\frac{q_{2k}}{m} \leq \frac{q_{2k}}{q_{j}} < \frac{q_{2k}}{q_{j-1}}\cdot \frac{q_{j-1}}{q_{j}}< \frac{1}{F_{j-1-2k}} \cdot \frac{1}{a_{j}} \to 0.\]
Thus
    \[
    \frac{q_{2k}}{m} + \frac{1}{2q_{2k}} \to 0.
    \]
We have shown that there is exactly one $a_i$ that is unbounded, and that $i$ is even.

We have $\frac rm = [0;a_1,\dots,a_{2k},\dots,a_n]$, with all of the $a_i$ fixed except $a_{2k}$, and $a_{2k}\to\infty$. Now
    \begin{align*}
    \lim H(r/m) &= \lim_{a_{2k}\to\infty} \; \min_{0\le j <n/2} q_{2j} \frac{r+1}m - p_{2j} \\
        &=\lim_{a_{2k}\to\infty} \; \min_{0\le j <n/2} \left( \frac{q_{2j}p_n-p_{2j}q_n+q_{2j}}{q_n} \right) \\
        &= \min_{0\le j <n/2} \; \lim_{a_{2k}\to\infty}
            \left(\frac{E[a_{2j+2},\dots,a_n]+E[a_1,\dots,a_{2j}]}{E[a_1,\dots,a_n]}\right)
    \end{align*}
Using the general identity (for $s\le \ell \le t$)
    \begin{multline*}
    E[a_s,\dots,a_t]=a_\ell E[a_s,\dots,a_{\ell-1}]E[a_{\ell+1},\dots,a_t]+\\
                        E[a_s,\dots,a_{\ell-2}]E[a_\ell+1,\dots,a_t]+E[a_s,\dots,a_{\ell-1}]E[a_{\ell+2},\dots,a_t]
    \end{multline*}
with $\ell=2k$, we can evaluate the limit as $a_{2k}\to\infty$. We arrive at
    \begin{align*}
    \beta = \lim H(\tfrac rm) &= \min\bigg\{
        \min_{0\le j <k}
            \frac{E[a_{2j+2},\dots,a_{2k-1}]E[a_{2k+1},\dots,a_n]}{E[a_1,\dots,a_{2k-1}]E[a_{2k+1},\dots,a_n]},\\
            & \hspace{5cm} \min_{k\le j < n/2}
            \frac{E[a_1,\dots,a_{2k-1}]E[a_{2k+1},\dots,a_{2j}}{E[a_1,\dots,a_{2k-1}]E[a_{2k+1},\dots,a_n]}
            \bigg\}\\
            &= \min\bigg\{
                \min_{0\le j <k}
                    \frac{E[a_{2j+2},\dots,a_{2k-1}]}{E[a_1,\dots,a_{2k-1}]},
            \min_{k\le j < n/2}
                \frac{E[a_{2k+1},\dots,a_{2j}}{E[a_{2k+1},\dots,a_n]}
            \bigg\}\\
            &= \min \bigg\{
                \frac{1}{E[a_1,\dots,a_{2k-1}]},
                \frac{1}{E[a_{2k+1},\dots,a_n]} \bigg\}  .
    \end{align*}
In either case, the numerator of $\beta$ is 1, and the proof of Theorem~\ref{MainTheorem} is concluded.

We note that we have actually proved (with a small bit of additional algebra) a quantitative version of the Theorem.
\begin{thm}
Let $(r_i,m_i)$ be a sequence of pairs of positive integers with $\gcd(r_i,m_i)=1$, $m_i\to\infty$ and $\limsup H(r_i/m_i)>0$. Then there is a pair of relatively prime positive integers $a,b$, with $a\le b$, a positive integer $c$, and an increasing sequence $i_1,i_2,\dots$ with
    \[r_{i_j} = \frac{am_{i_j}-c}{b} \qquad \text{ and }\qquad H\left(\frac{r_{i_j}}{m_{i_j}}\right) \to \frac{1}{\max\{c,b\}}.\]
\end{thm}

\end{document}